\documentclass[11pt,reqno]{amsart}
\usepackage{mathrsfs}

\usepackage{hyperref}

\usepackage{graphicx}
\usepackage{xcolor}
\usepackage{enumerate}
\usepackage{bm}

\newtheorem{theorem}{Theorem}[section]
\newtheorem{lemma}[theorem]{Lemma}
\newtheorem{proposition}[theorem]{Proposition}

 \theoremstyle{definition}

\theoremstyle{remark}

\numberwithin{equation}{section}

\begin{document}
\setlength{\baselineskip}{1.2\baselineskip}

\title[Interior gradient estimates in hyperbolic space]
{Interior gradient estimates for prescribed curvature equations in hyperbolic space}
%\date{2020.07.26}
%\thanks{}

\author{Zhenan Sui}

\address{Institute for Advanced Study in Mathematics of HIT, Harbin Institute of Technology, Harbin, China}
\email{sui.4@osu.edu}

\author{Wei Sun}

\address{Institue of Mathematical Sciences, ShanghaiTech University, Shanghai, China}
\email{sunwei@shanghaitech.edu.cn}

\begin{abstract}
In this paper, we study the interior gradient estimates for admissible solutions to prescribed curvature equations in hyperbolic space.
\end{abstract}

\maketitle

\section{Introduction}

In this note, we shall derive the interior gradient estimates for general prescribed curvature equations in hyperbolic space, which extends the recent work by Weng \cite{Weng2019} for $k$-th Weingarten curvature equations and by Mei-Zhou \cite{Mei-Zhou21} for curvature quotient equations. The techniques in this paper are based on our recent result \cite{SuiSun20} in Euclidean space. We shall elaborate each step towards the estimates to indicate the differences from Euclidean space.

Interior gradient estimate is a crucial step towards finding solutions with Lipschitz regularity on noncompact domains or domains with less regular boundary data. Around 1980s, Caffarelli-Nirenberg-Spruck proposed a question to find such estimate for prescribed curvature equations in their seminal paper ~\cite{CNS5}. Since then, interior gradient estimates have become an important topic in fully nonlinear elliptic equations.
Early influential work owes to Korevaar \cite{Korevaar87}, who used normal perturbation method to conduct computations on hypersurfaces subject to $k$-th Weingarten curvature equations. Korevaar's method was later generalized by Li \cite{LiYanyan91}, Trudinger \cite{Trudinger90} and  Guan-Spruck \cite{Guan-Spruck91} to general prescribed curvature equations as well as their parabolic version (see also \cite{Chou-Wang01,Chen15} for Hessian type equations).  In 1998,  Wang \cite{WangXJ98} gave a simplified PDE derivation of such estimate for $k$-th Weingarten curvature equations. Based on Wang's approach, Chen-Xu-Zhang \cite{ChenXuZhang17} derived the estimate for curvature quotient equations. In our paper \cite{SuiSun20}, we generalized this PDE method to a large class of prescribed curvature equations, including all previous equations and covering the degenerate cases. Our assumptions are motivated by the works\cite{LiYanyan91,Trudinger90}, and seem to be weaker. Indeed, our approach only utilizes fundamental structure conditions, combined with basic elements in convex analysis and linear algebra, from which one can see the essence of the derivation of such estimate.
In this note, we shall extend this approach to hyperbolic space by applying the crucial inequalities derived in \cite{SuiSun20}.

For hyperbolic space, we shall adopt the half space model
\[\mathbb{H}^{n+1} = \{ (x, x_{n+1}) \in \mathbb{R}^{n+1} \big\vert x_{n+1} > 0\}, \]
equipped with the metric
\[ d s^2 = x_{n+1}^{- 2} \sum_{i = 1}^{n+1} d x_i^2. \]

For a $C^2$ positive function $u$ defined on some domain in $\mathbb{R}^n$, we denote by $\bm{\lambda} (A) = (\lambda_1, \cdots, \lambda_n)$ the eigenvalue set of matrix $A$ with entries
\begin{equation}
	a_{ij} =  \frac{1}{w} \sum_{k,l} \gamma^{ik} (\delta_{kl} + u_k u_l + u u_{kl}) \gamma^{lj},
\end{equation}
%\begin{equation}
%	a_{ij} = \frac{1}{w} \left[\delta_{ij} + u \left( u_{ij} - \sum_l \frac{u_l u_j u_{li}}{w (w + 1)} - \sum_l \frac{u_l u_i u_{lj}}{w (w + 1)} + \sum_{k,l} \frac{u_i u_j u_k u_l u_{kl} }{w^2 (w + 1)^2} \right)\right],
%\end{equation}
where $w = \sqrt{1 + |D u|^2}$ and
\begin{equation*}
	\gamma^{ik} = \delta_{ik} - \frac{u_i u_k}{w (w + 1)} .
\end{equation*}
Following the literature \cite{CNS5} and in view of the calculations in \cite{GSS09}, $\lambda_1, \cdots, \lambda_n$ are called the principal curvatures of the vertical graph of $u$ in $\mathbb{H}^{n + 1}$ with respect to the upward normal vector field.

In this paper, we shall study prescribed curvature equations
\begin{equation}
	\label{main-equation}
		F(A) = \psi (x, u , \bm{\nu}) ,
\end{equation}
where 
\begin{equation*}
	F(A) = f (\bm{\lambda} (A)) , \quad\text{ and } \quad \bm{\nu} = \left( - \frac{D u}{w} , \frac{1}{w}\right) .
\end{equation*}
The right-side term $ \psi (x, u , \bm{\nu})$ is a prescribed $C^1$ function defined on some domain contained in
$\partial_{\infty} \mathbb{H}^{n+1} := \mathbb{R}^n \times \{0\} \cong \mathbb{R}^n $, 
and $f$ is a symmetric $C^1$ function defined on some open symmetric convex cone $\Gamma \subset \mathbb{R}^n$ with vertex at the origin containing the positive cone $\Gamma_n$.
We call a $C^2$ function $u$ {\em admissible} if $\bm{\lambda} (A) \in \bar \Gamma$. When $\bm{\lambda} (A) \in \partial \Gamma$, $f(\bm{\lambda} (A))$ should be interpreted as
$$
f(\bm{\lambda} (A)) = \limsup\limits_{\substack{ \bm{\lambda} \in \Gamma \\ \bm{\lambda} \rightarrow \bm{\lambda} (A) }} f(\bm{\lambda}). 
$$

As in the paper~\cite{SuiSun20}, we suppose that $f$ satisfies the structure conditions in $\Gamma$:
% as below:
\begin{enumerate}[(i) ]
\item $\frac{\partial f}{\partial\lambda_i} \geq 0$ for all $i$;
\item $f$ is concave;
\item  $\exists\, \frac{1}{2} > \delta > 0$, $\frac{\partial f}{\partial \lambda_j} > 2 \delta \sum_i \frac{\partial f}{\partial \lambda_i}$ when $\lambda_j < 0$;
\item $ f(t\bm{\lambda}) $ is increasing as a function of $t$ when %any $ \bm{\lambda} \in \Gamma$ and
$t > 0$;
\item $\lim_{t \rightarrow +\infty} f (t, \cdots ,t) > \sup_{B_r(0)}   \psi $.
\end{enumerate}
For the case that $\bm{\lambda} (A) \in \partial \Gamma$, 
instead of equation~\eqref{main-equation} we shall study, 
\begin{equation}
\label{main-equation-2}
	\tilde F(A) = h \circ \psi ,
\end{equation}
where
\begin{equation}
\label{equiv-equation}
	\tilde  F(A) :=  \limsup_{\substack{\bm{\lambda} \in \Gamma \\ \bm{\lambda} \to \bm{\lambda}(A)}}  g(\bm{\lambda}).
\end{equation}
In Equation~\eqref{main-equation-2}, function $g$ is assumed to be  continuously differentiable  in $\bar\Gamma \backslash (\partial \Gamma_n \cap \partial \Gamma)$, while function $h$ is assumed to be  strictly increasing and continuously differentiable on $\mathbb{R}$. The functions satisfies
\begin{equation}
\label{condition-g2}
	g(\bm{\lambda}) = h \circ f(\bm{\lambda}) \qquad \text{in } \Gamma .
\end{equation}
and
\begin{equation}
\label{condition-g1}
	\sum_i \frac{\partial g}{\partial \lambda_i} > 0, \qquad \text{when } \bm{\lambda} \in \partial \Gamma\backslash \partial \Gamma_n
\end{equation}
These assumptions will help us to deal with the boundary case $\bm{\lambda} (A) \in \partial\Gamma$. It comes that   $\bm{\lambda} (A) \in \partial \Gamma$ for many equations when $\psi \equiv f(\bm{0})$. 
If it always holds true that $\bm{\lambda}(A) \in \Gamma$, condition \eqref{condition-g1}--\eqref{condition-g2} is not in need as \eqref{main-equation} and \eqref{main-equation-2} are equivalent. 
A typical example of $f$  is the $k$th elementary symmetric function
\begin{equation*}
\sigma_k (\bm{\lambda}) :=  \sum_{1 \leq i_1 < \cdots < i_k \leq n} \lambda_{i_1} \cdots \lambda_{i_k}
\end{equation*}
defined on the $k$th G\r{a}rding's cone
\( \Gamma_k := \big\{\bm{\lambda} \in \mathbb{R}^n | \sigma_i (\bm{\lambda}) > 0, \, i = 1, \cdots, k \big\} \).
%When $f = \sigma_k^{\frac{1}{k}}$, we have to choose $g = \sigma_k$ and $h(t) = t^k$.
Moreover, the curvature quotients
\begin{equation*}
\Big(\frac{\sigma_k}{\sigma_l}\Big)^\frac{1}{k - l} (\bm{\lambda}), \quad 1 \leq l < k \leq n
\end{equation*}
also satisfy the structure conditions.

These structure conditions are well discussed in \cite{SuiSun20}, and we will only give a brief statement here. 
Condition (i)-(ii) are widely accepted in the study of fully nonlinear elliptic equations since Caffarelli-Nirenberg-Spruck~\cite{CNS3}.  
Condition (iii) was first proposed by Li \cite{LiYanyan91}, which was verified in \cite{LiYanyan91} for $f = \sigma_k^{\frac{1}{k}}$ and in \cite{Trudinger90} for $f = \left(\frac{\sigma_k}{\sigma_l} \right)^{\frac{1}{k - l}}$. It is crucial for interior gradient estimates in some geometric settings, although it is not necessary for fully nonlinear elliptic equations on Riemannian manifolds $(M,g)$~\cite{Sun20}, 
$$
F(\chi +  \nabla^2 u) = f(\bm{\lambda} (\chi + \nabla^2 u)) = \psi (u , du) .
$$
Condition (iv) implies that
\begin{equation}
	\sum_{i,j} \frac{\partial F}{\partial a_{ij}} a_{ij} \geq 0, \quad \text{whenever } \bm{\lambda} (A) \in \Gamma .
\end{equation}
From condition (v) and the simple fact that $\psi$ is locally defined on a compact set when $u$ is locally bounded, we know that there are constants $L > 0$, $\epsilon > 0$ depending on $\sup_{B_r (\bm{0})} \psi$ such that
\begin{equation}
	f (L, L, \cdots , L) > \sup_{B_r (\bm{0})} \psi + \epsilon .
\end{equation}
and hence by concavity,
\begin{equation}
	\sum_i \frac{\partial F}{\partial a_{ii}} > \frac{\epsilon}{L} .
\end{equation}
However, the coefficients are dependent on $\sup_{B_r(\bm{0})} \psi$.
A particular case is that $f$ is homogeneous of degree 1, which satisfies condition (iv) and (v). In this case,
\begin{equation}
	\sum_{i,j} \frac{\partial F}{\partial a_{ij}} a_{ij} = \psi ,
\end{equation}
and
\begin{equation}
	\sum_{i} \frac{\partial F}{\partial a_{ii}} = \psi + \sum_{i,j} \frac{\partial F}{\partial a_{ij}} (\delta_{ij} - a_{ij}) \geq f (\bm{1}) .
\end{equation}
Moreover, the coefficients are independent from $\sup_{B_r(\bm{0})} \psi$, and consequently the argument will be simpler.

The main result in this paper is the following interior gradient estimate for general prescribed curvature equation~\eqref{main-equation}.
\begin{theorem}
	Suppose that when $u > 0$,
\begin{equation}
	\label{condition-gradient}
	\psi_u - \frac{\psi - f(\bm{0})}{u}\geq 0 ,
	\end{equation}
where
 $f(\bm{0}) := \limsup_{\bm{\lambda} \to \bm{0}} f(\bm{\lambda}) $. Let $u \in C^3 (B_r (0))$ be an admissible solution to equation~\eqref{main-equation}. Then
 \begin{equation}
 	|Du (0)| \leq \exp \left( C \left(M + \frac{M^2}{r^2} + 1\right)\right)
 \end{equation}
where $M : = \sup_{B_r (\bm{0})} u$.

\end{theorem}

\medskip
\section{The gradient estimate}

In this section, we shall prove the interior gradient estimate.

%Now we begin to prove the interior gradient estimate. 
Let $u \in C^3 (B_r (0))$ be an admissible solution to equation~\eqref{main-equation}. Let us consider the test function
\begin{equation}
	G(x) := \ln \ln w + \varphi (u) + \ln \rho(x) ,
\end{equation}
where
\begin{equation}
	\varphi (u) = \ln u
	\end{equation}
and
\begin{equation}
	\rho(x) = (r^2 - |x|^2)^+ .
	\end{equation}

Suppose that $G(x)$ attains its maximal value at some point $p \in B_r (0)$. 
%If $u(p) = 0$, then it has to be that $D u (p) = \bm{0}$, which is a contradiction.
As in Sui-Sun~\cite{SuiSun20}, we are able to  choose a set of coordinates $(x_1, \cdots , x_n)$ such that
$	u_1 (p) = |D u (p)| $
and
$[u_{ij} (p)]_{i \geq 2, j\geq 2}$ is diagonal. Moreover, it is reasonable to require that $D^2 u(p)$ is in the form,
\begin{equation}
	\label{assumption-21}
	\left[
	\begin{array}{c|ccc|ccc}
		u_{11} & u_{12} & \cdots &u_{1m}   & 0   		&  \cdots   & 0		\\
		\hline
		u_{21} & u_{22} &            &  		  &  			&	   		   &                 \\
		\vdots &           & \ddots &            &     			&        	  &               \\
		u_{m1} &           &           &u_{mm} &  			  &			     &              \\
		\hline
		0 		  &			 &			  &		       & u_{m+1,m+1}	  &			     &			\\
		\vdots &           & 			&           &  				& 	\ddots &              \\
		0 &           &  			&			&            	&  			   &  u_{nn}
	\end{array}
	\right]
\end{equation}
with $u_{22} < u_{33} < \cdots < u_{mm}$ and nonzero entries $u_{21} , \cdots , u_{m1}$. It is easy to see that matrix $A$ at point $p$ is
\begin{equation}
	\left[
	\begin{array}{c|ccc|ccc}
		\dfrac{w^2 + u u_{11}}{w^3}  & \dfrac{u u_{12}}{w^2}  & \cdots & \dfrac{u u_{1m}}{w^2}    & 0   		&  \cdots   & 0		\\
		\hline
		\dfrac{ u u_{21}}{w^2} & \dfrac{1 + u u_{22}}{w}  &            &  		  &  			&	   		   &                 \\
		\vdots &           & \ddots &            &     			&        	  &               \\
		\dfrac{ u u_{m1}}{w^2} &           &           &\dfrac{ 1 + u u_{mm} }{w} &  			  &			     &              \\
		\hline
		0 		  &			 &			  &		       & \dfrac{1 + u u_{m+1,m+1}	}{w}  &			     &			\\
		\vdots &           & 			&           &  				& 	\ddots &              \\
		0 &           &  			&			&            	&  			   & \dfrac{1 + u u_{nn}}{w}
	\end{array}
	\right] .
\end{equation}
Thus, $a_{22} < a_{33} < \cdots < a_{mm}$ and $a_{21}, \cdots , a_{m1}$ are nonzero. Moreover
\begin{equation}
\begin{aligned}
 F^{11} &:= \frac{\partial F}{\partial u_{11}} = \frac{u}{w^3} \frac{\partial F}{\partial a_{11}} ,  &&\quad \\
 F^{1i} &:= \frac{\partial F}{\partial u_{1i}} = \frac{u}{w^2} \frac{\partial F}{\partial a_{1i}} ,   &&\quad 2 \leq i \leq n \\
F^{ij} &:= \frac{\partial F}{\partial u_{ij}} = \frac{u}{w} \frac{\partial F}{\partial a_{ij}} ,  &&\quad 2 \leq i,j \leq n .
\end{aligned}
\end{equation}

There are two cases: (1) $\bm{\lambda} (A) \in \Gamma$ and $\bm{\lambda} (A) \in \partial\Gamma$.

\medskip
\subsection{Case 1.} $\bm{\lambda} (A) \in \Gamma$ at point $p$. Let us first recall some inequalities from Sui-Sun~\cite{SuiSun20}.

By symmetry and concavity of $F$, it  is easy to obtain the following inequality.
\begin{lemma}
	\label{lemma-1}
	Under the particular set of coordinates,
	$$
	\frac{\partial F}{\partial a_{1i}} a_{1i} \leq 0 , \quad \forall \, i \geq 2.
	$$
\end{lemma}

%The coefficient in Lemma~\ref{lemma-2} might be different than that in condition (iii). 
Generally, some kind of boundary condition is imposed to guarantee that the solution is admissible in the study of fully nonlinear elliptic equations, e.g.
$$
\limsup_{\substack{\bm{\lambda} \to \Gamma \\ \bm{\lambda } \to \bm{\lambda_0}}}
 f(\bm{\lambda}) < \inf \psi, \qquad \forall \bm{\lambda_0} \in \partial \Gamma,
 $$
or
$$
f > 0 \quad \text{ in } \Gamma, \qquad f = 0 \quad \text{on } \Gamma .
$$
We are able to avoid this kind of assumption in the gradient estimate. By an argument via dual cone, we can prove a subtle property of convex cone $\Gamma$.
\begin{proposition}
\label{proposition-1}
%Let $\Gamma$  be a convex symmetric cone with vertex at the origin and containing $(1,1,\cdots , 1)$. Then 
There exists a constant $\delta' > 0$ such that if $\bm{\lambda} = (\lambda_1 ,\cdots ,\lambda_n) \in \bar\Gamma$ with $\lambda_1 < 0$, we have
$$
(0, \lambda_2 + \delta' \lambda_1 , \cdots , \lambda_n + \delta' \lambda_1) \in \Gamma .
$$
\end{proposition}

Without loss of generality, we may assume that $\delta \leq \frac{\delta'}{1 + \delta'}$. By the concavity, condition~(iii) and Proposition~\ref{proposition-1}, we are able to prove the following inequality.
\begin{lemma}
	\label{lemma-2}
	Under the particular set of coordinates, if $a_{11} < 0$, then we have
	\begin{equation}
		\frac{\partial F}{\partial a_{11}} a^2_{11} + 2 \sum_{i \geq 2} \frac{\partial F}{\partial a_{1i}} a_{1i} a_{11}
		\geq \delta a^2_{11} \sum_i \frac{\partial F}{\partial a_{ii}} .
	\end{equation}
\end{lemma}

Chou-Wang~\cite{Chou-Wang01} adopted an inequality for $\sigma_k$, while Chen~\cite{Chen15} found a similar inequality for $\frac{\sigma_k}{\sigma_l}$. Both of their arguments rely on the special properties of $\sigma_k$, especially the fact that $\sigma_k (A)$ is actually a polynomial of entries of $A$. We were able to prove the following inequality by calculation of eigenvalues.
\begin{lemma}
	\label{lemma-3}
	Under the particular set of coordinates,
	\begin{equation}
		\sum_{i\geq 2} \frac{\partial F}{\partial a_{1i}} a_{1i} a_{ii} \geq - \left(\frac{a^2_{11}}{4} + \sum_{i\geq 2} a^2_{1i}\right) \frac{\partial F}{\partial a_{11}} .
	\end{equation}
\end{lemma}

Now we can apply the  maximum principle. Differentiating test function $G$ twice, we have
\begin{equation}
	\label{G-1}
%	0 = 
	G_i = \frac{w_i}{w \ln w} + \varphi' u_i + \frac{\rho_i}{\rho} = \frac{u_1 u_{1i}}{w^2 \ln w} + \varphi' u_i + \frac{\rho_i}{\rho} ,
\end{equation}
and
\begin{equation}
	\label{G-2}
	\begin{aligned}
		G_{ij}
%		&= \frac{w_{ij}}{w \ln w} - \frac{(1 + \ln w) w_i w_j}{w^2 \ln^2 w} + \varphi' u_{ij} + \varphi'' u_i u_j + \frac{\rho_{ij}}{\rho} - \frac{\rho_i \rho_j}{\rho^2} \\
%		&=  \frac{ u_1 u_{1ij}}{ w^2 \ln w}+  \frac{\sum_k u_{ki} u_{kj} }{ w^2 \ln w} - \frac{  2 u^2_1 u_{1i}  u_{1j}}{w^4 \ln w} - \frac{ u^2_1 u_{1i} u_{1j}}{w^4 \ln^2 w}  \\
%		&\qquad + \varphi' u_{ij} + \varphi'' u_i u_j + (\varphi')^2 u_i u_j + \frac{\rho_{ij}}{\rho}  - \frac{u^2_1 u_{1i} u_{1j}}{w^4 \ln^2 w} +  \frac{\varphi'}{\rho} (u_i \rho_j + u_j \rho_i) \\
		&=  \frac{ u_1 u_{1ij}}{ w^2 \ln w}+  \frac{\sum_k u_{ki} u_{kj} }{ w^2 \ln w} - \frac{  2 u^2_1 u_{1i}  u_{1j}}{w^4 \ln w} - \frac{ 2 u^2_1 u_{1i} u_{1j}}{w^4 \ln^2 w}  \\
		&\qquad + \varphi' u_{ij} + \frac{\rho_{ij}}{\rho}  +  \frac{\varphi'}{\rho} (u_i \rho_j + u_j \rho_i) .
	\end{aligned}
\end{equation}
Rewriting \eqref{G-2}, we obtain
\begin{equation}
	\label{G-2-1}
	\begin{aligned}
		G_{ij}
		&= \frac{u_1 u_{1ij}}{w^2 \ln w} + \frac{1}{u^2 \ln w}   \left(w^2 - 2u^2_1 - \frac{2 u^2_1}{\ln w}\right) \gamma_{ii} \gamma_{jj}  a_{1i} a_{1j} \\
		&\qquad + \frac{1}{u^2 w^4\ln w}   \left(w^2 - 2u^2_1 - \frac{2 u^2_1}{\ln w}\right)    (\delta_{1i} + u_1 u_i)   (\delta_{1j} + u_1 u_j) \\
		&\qquad - \frac{1}{u^2 w^2 \ln w}   \left(w^2 - 2u^2_1 - \frac{2 u^2_1}{\ln w}\right) \gamma_{ii}   a_{1i} (\delta_{1j} + u_1 u_j) \\
		&\qquad - \frac{1}{u^2 w^2 \ln w}   \left(w^2 - 2u^2_1 - \frac{2 u^2_1}{\ln w}\right)  \gamma_{jj}   (\delta_{1i} + u_1 u_i) a_{1j}  \\
		&\qquad + \frac{1}{u^2 \ln w}  \sum_{k\geq 2}  \gamma_{ii} \gamma_{jj}  a_{ki} a_{kj}   + \frac{1}{u^2 w^2 \ln w}  \sum_{k\geq 2}   \delta_{ki}   \delta_{kj}  \\
		&\qquad	 - \frac{1}{u^2 w \ln w}  \sum_{k\geq 2}  \gamma_{ii}  a_{ki}     \delta_{kj}  - \frac{1}{u^2 w \ln w}  \sum_{k\geq 2}   \gamma_{jj}   \delta_{ki} a_{kj}
		\\
		&\qquad + \varphi' \frac{w}{u} \gamma_{ii} \gamma_{jj} a_{ij}  -  \frac{\varphi'}{u}  (\delta_{ij} + u_i u_j)  - \frac{2 \delta_{ij}}{\rho} + \frac{\varphi'}{\rho} (u_i \rho_j + u_j \rho_i) .
	\end{aligned}
\end{equation}
Differentiating equation \eqref{main-equation}, we can see that at point $p$,
\begin{equation}
	\label{F-1}
 	\begin{aligned}
 		\sum_{i,j} F^{ij} u_{ij1}
 		&= \partial_1 \psi - \left(\varphi' u_1 + \frac{\rho_1}{\rho}\right)  \ln w  \sum_{i,j} \frac{\partial F}{\partial a_{ij}}  a_{ij}   - \frac{ u_1 }{u }  \sum_{i,j} \frac{\partial F}{\partial a_{ij}}  a_{ij} +  \frac{2 u_1 }{u w} \frac{\partial F}{\partial a_{11}}   \\
 		&\qquad  + \frac{2 u_1 w }{u}  \frac{\partial F}{\partial a_{11}}  a^2_{11}  + \frac{2 u_1 w (2 w + 1)}{u(w + 1)} \sum_{i \geq 2} \frac{\partial F}{\partial a_{i1}}  a_{11}  a_{i1} \\
 		&\qquad  + \frac{2 u_1 w}{u  (w + 1)}   \frac{\partial F}{\partial a_{11}}   \sum_{i \geq 2}   a_{i1}  a_{i1}  + \frac{2 u_1 w^2}{u(w + 1)} \sum_{i\geq 2,j\geq 2} \frac{\partial F}{\partial a_{ij}}  a_{1j}  a_{i1} \\
 		&\qquad + \frac{2 u_1 w}{u  (w + 1)}  \sum_{j \geq 2} \frac{\partial F}{\partial a_{1j}}     a_{jj}  a_{j1} -  \frac{4 u_1  }{u }  \sum_{j} \frac{\partial F}{\partial a_{1j}}  a_{1j} + \frac{u_1}{u w} \sum_{i} \frac{\partial F}{\partial a_{ii}} .
 	\end{aligned}
 \end{equation}
From \eqref{G-1}, at point $p$,
\begin{equation}
	\label{G-1-1}
		\frac{u_1 u_{11}}{w^2 \ln w} = - \varphi' u_1 + \frac{2 x_1}{\rho} ,
	\end{equation}
and
\begin{equation}
	\label{G-1-2}
		\frac{u_1 u_{1i}}{w^2 \ln w} = \frac{2 x_i}{\rho} , \qquad \text{for } i \geq 2 .
	\end{equation}
We may assume that at point $p$,
\begin{equation}
\label{G-1-3}
	\rho	\ln |D u| \geq 1000 r^2+ 4 Mr .
\end{equation}
%and hence
%\begin{equation}
%	\ln w \geq \ln |D u| \geq C + \frac{4 M}{r} .
%\end{equation}
Otherwise, the proof is finished. By  \eqref{G-1-3} and \eqref{G-1-1},
\begin{equation}
	u_{11}= \frac{w^2 \ln w}{\rho u_1} \left( -\varphi' \rho u_1 + {2 x_1}\right) \leq - \frac{\varphi'}{2} w^2 \ln w
	= - \frac{1}{2 u} w^2 \ln w ,
\end{equation}
and hence
\begin{equation}
\label{G-1-4}
%\begin{aligned}
	a_{11}
%	&
	= \frac{w^2 + u u_{11}}{w^3}
%	\\
%	&
	\leq \frac{1}{w} \left(1 - \frac{u \varphi' \ln w}{2}\right) 
%	\\
%	&=
%		\frac{1}{w} \left(1 - \frac{\ln w}{2}\right)\\
%	&
	\leq - \frac{\ln w}{4 w} .
%\end{aligned}
\end{equation}

From \eqref{G-2}, at point $p$,
\begin{equation}
	\label{max-1}
	\begin{aligned}
		0
		&\geq 
%		\sum_{i,j} F^{ij} G_{ij} \\
%		&= 
		\frac{u_1 }{w^2 \ln w} \sum_{i,j} F^{ij} u_{1ij} + \frac{1}{u w \ln w} \sum_{\substack{i,j,k }}  \frac{\partial F}{\partial a_{ij}}     a_{ki} a_{kj} \\
		&\qquad + \frac{1}{u w \ln w}   \left( - u^2_1 - \frac{2 u^2_1}{\ln w}\right) \sum_{i,j}  \frac{\partial F}{\partial a_{ij}}   a_{1i} a_{1j}   \\
		&\qquad + \frac{1}{u w^3\ln w}   \left( - u^2_1 - \frac{2 u^2_1}{\ln w}\right)  \frac{\partial F}{\partial a_{11}}   + \frac{1}{u w^3 \ln w} \sum_{\substack{ i }} \frac{\partial F}{\partial a_{ii}}        \\
		&\qquad - \frac{2}{u w^2 \ln w}   \left( - u^2_1 - \frac{2 u^2_1}{\ln w}\right) \sum_{i}  \frac{\partial F}{\partial a_{i1}}    a_{1i}     - \frac{2}{u w^2 \ln w}  \sum_{\substack{i,j}}  \frac{\partial F}{\partial a_{ij}}     a_{ji} \\
		&\qquad + \frac{1}{u}  \sum_{i,j}  \frac{\partial F}{\partial a_{ij}}  a_{ij}   -  \frac{1}{u w} \sum_{i}  \frac{\partial F}{\partial a_{ii}} (\gamma^{ii})^2    -  \frac{u^2_1}{u w^3}   \frac{\partial F}{\partial a_{11}}  \\
		&\qquad   - \frac{2u}{\rho w} \sum_{i} \frac{\partial F}{\partial a_{ii}} (\gamma^{ii} )^2 + \frac{2   }{\rho w} \sum_{i}  \frac{\partial F}{\partial a_{1i}} \gamma^{ii}  u_1 \rho_i  .
		\end{aligned}
	\end{equation}
Substituting \eqref{F-1} and \eqref{G-1-2} into \eqref{max-1} and applying Lemma~\ref{lemma-1}, \ref{lemma-2} and \ref{lemma-3},
\begin{equation}
\label{inequality-6}
	\begin{aligned}
	0
	&\geq  \frac{u_1 \psi_{x_1} }{w^2 \ln w} + \frac{u^2_1 \psi_u }{w^2 \ln w} - \frac{u_1 \psi_{\nu_1} u_{11}}{w^5 \ln w} - \frac{u_1 \sum_{i\geq 2} \psi_{\nu_i} u_{i1}}{w^3 \ln w}  - \frac{ \psi_{\nu_{n + 1}} u^2_1 u_{11}}{w^5 \ln w}  \\
	&\quad   -\frac{u_1 \rho_1}{ \rho w^2 } \sum_{i,j} \frac{\partial F}{\partial a_{ij}}  a_{ij}   - \frac{w^2 + 1}{u w^2  \ln w}  \sum_{\substack{i,j}}  \frac{\partial F}{\partial a_{ij}}     a_{ij}  + \frac{1}{u w^2}  \sum_{i,j}  \frac{\partial F}{\partial a_{ij}}  a_{ij}     \\
	&\quad + \frac{\delta w}{2 u  \ln w}  a^2_{11}  \sum_i \frac{\partial F}{\partial a_{ii}} -  \frac{ 1 }{u  w   \ln w}   \frac{\partial F}{\partial a_{11}}   \sum_{i \geq 2}   a^2_{i1}   \\
	&\quad     -  \frac{2}{u w} \sum_{i}  \frac{\partial F}{\partial a_{ii}}       - \frac{2u}{\rho w} \sum_{i} \frac{\partial F}{\partial a_{ii}}    + \frac{2  \rho_1  u_1 }{\rho w^3}  \frac{\partial F}{\partial a_{11}}  .
	\end{aligned}
\end{equation}

Note that
\begin{equation}
%	\begin{aligned}
	f(\bm{\eta}) - \psi \leq \sum_{i,j} \frac{\partial F}{\partial a_{ij}} \left(\eta_i \delta_{ij} - a_{ij}\right) = \sum_i \frac{\partial F}{\partial a_{ii}} \eta_i - \sum_{i,j} \frac{\partial F}{\partial a_{ij}} a_{ij} ,
%	\end{aligned}
	\end{equation}
and hence
\begin{equation}
	\sum_{i,j} \frac{\partial F}{\partial a_{ij}} a_{ij} \leq \sum_i \frac{\partial F}{\partial a_{ii}} \eta_i - f(\bm{\eta}) + \psi.
	\end{equation}
Choosing a sequence $\bm{\eta}_n \rightarrow \bm{0}$, it is possible to reach that
\begin{equation}
\label{concavity-1}
	\sum_{i,j} \frac{\partial F}{\partial a_{ij}} a_{ij} \leq \psi - f(\bm{0}) ,
	\end{equation}
and hence by assumption \eqref{condition-gradient} and inequality~\eqref{G-1-3}, we have
%\begin{equation}
%	\psi_u \geq \frac{1}{u} \sum_{i,j} \frac{\partial F}{\partial a_{ij}} a_{ij} .
%\end{equation}
\begin{equation}
\label{inequality-7}
\begin{aligned}
	\frac{u^2_1 \psi_u}{w^2 \ln w} - \frac{w^2 + 1}{u w^2 \ln w} \sum_{i,j} \frac{\partial F}{\partial a_{ij}} a_{ij} + \frac{1}{u w^2} \sum_{i,j} \frac{\partial F}{\partial a_{ij}} a_{ij} \\
	\geq \left( 1 - \frac{2}{  \ln w}  \right)\frac{1}{u w^2} \sum_{i,j} \frac{\partial F}{\partial a_{ij}} a_{ij} 
%	\geq \left( 1 - \frac{\rho}{  1000 r^2}  \right)\frac{1}{u w^2} \sum_{i,j} \frac{\partial F}{\partial a_{ij}} a_{ij} \\
	\geq 0.
\end{aligned}
\end{equation}
Using \eqref{G-1-1} and \eqref{G-1-2}, we have
\begin{equation}
\label{term-1}
\begin{aligned}
	 &\qquad \frac{u_1 \psi_{x_1} }{w^2 \ln w} - \frac{u_1 \psi_{\nu_1} u_{11}}{w^5 \ln w} - \frac{u_1 \sum_{i\geq 2} \psi_{\nu_i} u_{i1}}{w^3 \ln w}  - \frac{ \psi_{\nu_{n + 1}} u^2_1 u_{11}}{w^5 \ln w}   \\
 	&= \frac{u_1 \psi_{x_1} }{w^2 \ln w} - \frac{ \psi_{\nu_1} }{w^3} \left( - \frac{u_1}{u} + \frac{2 x_1}{\rho} \right) - \frac{ 2 \sum_{i\geq 2} \psi_{\nu_i} x_i}{\rho w}  - \frac{ \psi_{\nu_{n + 1}} u_1}{w^3} \left( - \frac{u_1}{u} + \frac{2 x_1}{\rho} \right) \\
	 	&\geq - \frac{|D \psi| u_1}{w^2 \ln w} - \frac{|D \psi|}{w^3} \left(\frac{u_1}{u} + \frac{2 r}{\rho}\right) - \frac{2 |D \psi| r}{\rho w} - \frac{|D \psi| u_1}{w^3} \left(\frac{u_1}{u} + \frac{2 r}{\rho}\right) \\
	 		&\geq  - \frac{|D \psi| }{w \ln w}  - \frac{6 |D \psi| r}{\rho w} - \frac{2 |D \psi| }{u w} \\
	 		&\geq  - \frac{|D \psi| }{w }  - \frac{6 |D \psi| r}{\rho w} - \frac{2 |D \psi| }{u w} .
\end{aligned}
\end{equation}
Also by \eqref{G-1-2},
\begin{equation}
\label{term-2}
\begin{aligned}
	- \frac{1}{u w \ln w} \frac{\partial F}{\partial a_{11}} \sum_{i\geq 2} a^2_{i1}
	&= - \frac{1}{u w \ln w} \frac{\partial F}{\partial a_{11}} \sum_{i\geq 2} \frac{u^2}{w^4}  u^2_{i1} \\
	&\geq - \frac{4 r^2 u \ln w}{\rho^2  u^2_1 w }  \frac{\partial F}{\partial a_{11}}\\
	&\geq - \frac{5 r^2 u \ln w}{\rho^2  w^3 }  \frac{\partial F}{\partial a_{11}}.
\end{aligned}
\end{equation}
Substituting \eqref{G-1-2},  \eqref{concavity-1} -- \eqref{term-2} into \eqref{inequality-6},
\begin{equation}
\label{inequality-2}
	\begin{aligned}
	0
	&\geq  - \frac{|D \psi| }{w }  - \frac{6 |D \psi| r}{\rho w} - \frac{2 |D \psi| }{u w} -\frac{2 r u_1 }{ \rho w^2 } \left(\psi - f(\bm{0})\right)\\
	&\quad    +   \frac{\delta \ln w}{32 u w}  \sum_i \frac{\partial F}{\partial a_{ii}}   - \frac{5 r^2 u \ln w}{\rho^2  w^3 }  \frac{\partial F}{\partial a_{11}}\\
	&\quad   -  \frac{2}{u w} \sum_{i}  \frac{\partial F}{\partial a_{ii}}  - \frac{2u}{\rho w} \sum_{i} \frac{\partial F}{\partial a_{ii}}  - \frac{4 r u_1 }{ \rho w^3}  \frac{\partial F}{\partial a_{11}}  .
  	\end{aligned}
\end{equation}
Multiplying \eqref{inequality-2} by $\rho u w$, we obtain
\begin{equation}
\label{inequality-3}
\begin{aligned}
	\frac{\delta \rho \ln w}{32 }  \sum_i \frac{\partial F}{\partial a_{ii}}
	&\leq   \rho  M |D \psi| + 6 M |D \psi| r + 2 \rho |D \psi| + 2 r M  \left(\psi - f(\bm{0})\right)\\
	&      + \frac{5 r^2 M^2 \ln w}{\rho  w^2 }  \frac{\partial F}{\partial a_{11}} +  2 \rho \sum_{i}  \frac{\partial F}{\partial a_{ii}}  + 2 M^2 \sum_{i} \frac{\partial F}{\partial a_{ii}}  + \frac{4 r M  }{ w}  \frac{\partial F}{\partial a_{11}}  .
\end{aligned}
\end{equation}
Dividing \eqref{inequality-3} by $\sum_i \frac{\partial F}{\partial a_{ii}}$,
\begin{equation}
\begin{aligned}
	 \frac{\delta \rho \ln w}{32 }
	&\leq   \frac{\rho L}{\epsilon} M |D \psi| + \frac{6 L M |D \psi| r}{\epsilon } + \frac{2 \rho L}{\epsilon} |D \psi| + \frac{2 L r M }{ \epsilon  } \left(\psi - f(\bm{0})\right)\\
	&\quad    + \frac{5 r^2 M^2 \ln w}{\rho  w^2 } +  2 \rho  + 2 M^2   + \frac{4 r M  }{ w}   .
\end{aligned}
\end{equation}
Therefore,
\begin{equation}
	\ln w (\bm{0}) \leq C \left( M + \frac{M^2}{r^2} + 1\right) .
\end{equation}

\medskip
\subsection{Case 2.} $\bm{\lambda}(A) \in \partial\Gamma$. In this case, we consider the equivalent equation~\eqref{equiv-equation}.

By approximation, we have
\begin{equation}
	h'(F(A)) (\psi - f(\bm{0}))\geq \sum_{i,j} \frac{\partial \tilde F}{\partial a_{ij}} a_{ij} \geq 0 ,
\end{equation}
and
\begin{equation}
	\sum_{i} \frac{\partial \tilde F}{\partial a_{ii}} \geq h' \circ F(A) \frac{\epsilon}{L} .
\end{equation}
Similarly, we also have
\begin{equation}
	\frac{\partial \tilde F}{\partial a_{1i}} a_{1i} \leq 0 ,
	\end{equation}
and
\begin{equation}
	\frac{\partial \tilde F}{\partial a_{11}} a^2_{11} + 2 \sum_{i\geq 2} \frac{\partial \tilde F}{\partial a_{1i}} a_{1i} a_{11} \geq \delta a^2_{11} \sum_i \frac{\partial \tilde F}{\partial a_{ii}} ,
	\end{equation}
and
\begin{equation}
	\sum_{i\geq 2} \frac{\partial \tilde F}{\partial a_{1i}} a_{1i} a_{ii} \geq - \left(\frac{a^2_{11}}{4} + \sum_{i\geq 2} a^2_{1i}\right) \frac{\partial \tilde F}{\partial a_{11}} .
	\end{equation}

As in case 1, at point $p$,
\begin{equation}
	\label{inequality-4}
	\begin{aligned}
		0
		&\geq  \frac{u_1 h'(\psi) \psi_{x_1} }{w^2 \ln w} - \frac{u_1 h'(\psi) \psi_{\nu_1} u_{11}}{w^5 \ln w} - \frac{u_1 h'(\psi) \sum_{i\geq 2} \psi_{\nu_i} u_{i1}}{w^3 \ln w}  - \frac{ h'(\psi) \psi_{\nu_{n + 1}} u^2_1 u_{11}}{w^5 \ln w}  \\
		&\quad   -\frac{u_1 \rho_1}{ \rho w^2 } \sum_{i,j} \frac{\partial \tilde F}{\partial a_{ij}}  a_{ij}  + \frac{\delta w}{2 u  \ln w}  a^2_{11}  \sum_i \frac{\partial \tilde F}{\partial a_{ii}} -  \frac{ 1 }{u  w   \ln w}   \frac{\partial \tilde  F}{\partial a_{11}}   \sum_{i \geq 2}   a^2_{i1}   \\
		&\quad   + \frac{1}{u w \ln w} \sum_{\substack{ i }} \frac{\partial \tilde F}{\partial a_{ii}}  -  \frac{1}{u w} \sum_{i}  \frac{\partial \tilde F}{\partial a_{ii}}   -  \frac{ u^2_1}{u w^3}   \frac{\partial \tilde F}{\partial a_{11}}     - \frac{2u}{\rho w} \sum_{i} \frac{\partial \tilde F}{\partial a_{ii}}  + \frac{2  \rho_1  u_1 }{ \rho w^3}  \frac{\partial \tilde F}{\partial a_{11}}  \\
%		&\geq - \frac{h'(\psi) |D \psi|}{w} \left( 1 + \frac{6r}{\rho} + \frac{2}{u}\right) - \frac{2 r u_1}{\rho w^2} h'(\psi)\left(\psi - f(\bm{0})\right) - \frac{5 r^2 u \ln w}{\rho^2 w^3} \frac{\partial \tilde F}{\partial a_{11}} \\
%		&\quad    + \frac{\delta w  a^2_{11}}{2 u  \ln w}   \sum_i \frac{\partial \tilde F}{\partial a_{ii}}    + \frac{1}{u w \ln w} \sum_{\substack{ i }} \frac{\partial \tilde F}{\partial a_{ii}}  -  \frac{1}{u w} \sum_{i}  \frac{\partial \tilde F}{\partial a_{ii}}   -  \frac{ 1}{u w}   \frac{\partial \tilde F}{\partial a_{11}}     - \frac{2u}{\rho w} \sum_{i} \frac{\partial \tilde F}{\partial a_{ii}}  + \frac{2  \rho_1   }{ \rho w^2}  \frac{\partial \tilde F}{\partial a_{11}}  \\
		&\geq - \frac{h'(\psi) |D \psi|}{w} \left( 1 + \frac{6r}{\rho} + \frac{2}{u}\right) - \frac{2 r }{\rho w} h'(\psi)\left(\psi - f(\bm{0})\right) - \frac{5 r^2 u \ln w}{\rho^2 w^3}  \sum_i \frac{\partial \tilde F}{\partial a_{ii}}   \\
		&\quad    + \frac{\delta \ln w  }{3 2 u  w}   \sum_i \frac{\partial \tilde F}{\partial a_{ii}}      -  \frac{2}{u w} \sum_{i}  \frac{\partial \tilde F}{\partial a_{ii}}       - \frac{2u}{\rho w} \sum_{i} \frac{\partial \tilde F}{\partial a_{ii}}  - \frac{4 r   }{ \rho w^2}   \sum_i \frac{\partial \tilde F}{\partial a_{ii}}    .
	\end{aligned}
\end{equation}
Multiplying \eqref{inequality-4} by $\rho u w$,
\begin{equation}
	\label{inequality-5}
	\begin{aligned}
		 \frac{\delta \rho \ln w  }{3 2 }   \sum_i \frac{\partial \tilde F}{\partial a_{ii}}
		&\leq
		 h'(\psi) |D \psi| r^2 M + 6 h'(\psi) |D \psi| M r + 2 h'(\psi) |D \psi| r^2 \\
		 &\quad + 2 r M h'(\psi)\left(\psi - f(\bm{0})\right) + \frac{5 r^2 M^2 \ln w}{\rho w^2}  \sum_i \frac{\partial \tilde F}{\partial a_{ii}}   \\
		&\quad     +   2 r^2 \sum_{i}  \frac{\partial \tilde F}{\partial a_{ii}}   +  2 M^2  \sum_{i} \frac{\partial \tilde F}{\partial a_{ii}}  + \frac{4 M r   }{ w}   \sum_i \frac{\partial \tilde F}{\partial a_{ii}}    .
	\end{aligned}
\end{equation}
If $h'(\psi) = 0$,
\begin{equation}
	\begin{aligned}
		\frac{\delta \rho \ln w  }{3 2 }
		&\leq  \frac{5 r^2 M^2 \ln w}{\rho w^2}     +   2 r^2   +  2 M^2     + \frac{4 M r   }{ w}     ;
	\end{aligned}
\end{equation}
if $h'(\psi) > 0$,
\begin{equation}
	\begin{aligned}
		\frac{\delta \rho \ln w}{32 }
		&\leq   \frac{ L}{\epsilon} r^2M |D \psi| + \frac{6 L M |D \psi| r}{\epsilon } + \frac{2  L}{\epsilon} |D \psi| r^2 + \frac{2 L r M }{ \epsilon  } \left(\psi - f(\bm{0})\right)\\
		&\quad    + \frac{5 r^2 M^2 \ln w}{\rho  w^2 } +  2 r^2  + 2 M^2   + \frac{4  M r }{ w}   .
	\end{aligned}
\end{equation}
Therefore,
\begin{equation}
	\ln w (\bm{0}) \leq C \left( M + \frac{M^2}{r^2} + 1\right) .
\end{equation}

%\newpage
\medskip
\noindent
{\bf Acknowledgements}\quad
The first author is supported by NSFC grant (No. 12001138).
The second author is supported by NSFC grant (No. 11501119) and a start-up grant from ShanghaiTech University (2018F0303-000-04).

\end{document}